\documentclass[12pt]{article}

\usepackage{amsmath}
\usepackage{amssymb}
\usepackage[mathscr]{eucal}
\usepackage{graphicx}
\usepackage{color}

\newtheorem{Theorem}{\sc Theorem}
\newtheorem{Definition}[Theorem]{\sc Definition}
\newtheorem{Proposition}[Theorem]{\sc Proposition}
\newtheorem{Lemma}[Theorem]{\sc Lemma}
\newtheorem{Corollary}[Theorem]{\sc Corollary}

\newtheorem{Example}[Theorem]{\sc Example}

\newcommand{\cP}{\mbox{{${\cal P}$}}}

\newcommand{\ve}{\mbox{{$\varepsilon$}}}
\newcommand{\R}{{\if mm {\rm I}\mkern -3mu{\rm R}\else \leavevmode
\hbox{I}\kern -.17em\hbox{R} \fi}}

\newcommand{\cT}{\mbox{{${\cal T}$}}}
\newcommand{\cC}{\mbox{{${\cal C}$}}}

\newcommand{\cF}{\mbox{{${\cal F}$}}}

\newcommand{\bu}{\mbox{\boldmath{$u$}}}
\newcommand{\bv}{\mbox{\boldmath{$v$}}}

\newcommand{\bx}{\mbox{\boldmath{$x$}}}

\newcommand{\fb}{\mbox{\boldmath{$f$}}}

\newcommand{\bvarepsilon}{\mbox{\boldmath{$\varepsilon$}}}

\newcommand{\bnu}{\mbox{\boldmath{$\nu$}}}

\newcommand{\bzero}{\mbox{\boldmath{$0$}}}

\def\sqr#1#2{{
    \vcenter{
         \vbox{\hrule height.#2pt
               \hbox{\vrule width.#2pt height#1pt \kern#1pt
                     \vrule width.#2pt
               }
               \hrule height.#2pt
         }
    }
}}

\def\real{\mathbb{R}}

\def\lista#1
{{ \itemindent 0.0cm \labelsep .2cm \leftmargin 0.8cm \rightmargin
0.0cm \labelwidth 0.6cm \topsep 0.0mm
\parsep 0.0mm
\itemsep 0.0mm
\begin{list}{}
{ \setlength{\leftmargin}{.8cm} \setlength{\rightmargin}{0.0cm}
\setlength{\parsep}{0.0mm} \setlength{\topsep}{.0mm}
\setlength{\parskip}{.0cm} \setlength{\itemsep}{.0cm} }
{#1}\end{list}} }

\newcounter{theorem}

\begin{document}



\title{\bf  A Convergence Criterion for Elliptic Variational Inequalities}

\vspace{22mm}
{\author{Claudia Gariboldi$^{1}$,\ Anna Ochal$^{2}$,
\ Mircea Sofonea$^{3}$\footnote{Corresponding author, E-mail : sofonea@univ-perp.fr}\\ and\\ Domingo A. Tarzia$^{4,5}$\\[5mm]
{\it \small $1$ Departamento de Matem\'atica, FCEFQyN}\\[-1mm]
{\it \small
	Universidad Nacional de Rio Cuarto}
\\[-1mm]
{\it\small Ruta 36 km 601, 5800 Rio Cuarto, Argentina}		\\[5mm]	
{\it \small $2$ Chair in Optimization and Control}\\[-1mm]
{\it \small
	Jagiellonian University in Krakow}
\\[-1mm]
{\it\small Ul. Lojasiewicza 6, 30348 Krakow, Poland}		\\[5mm]	
{\it \small $3$ Laboratoire de Math\'ematiques et Physique}\\[-1mm]
{\it \small
	University of Perpignan Via Domitia}
\\[-1mm]{\it\small 52 Avenue Paul Alduy, 66860 Perpignan, France}		\\[6mm]		
{\it\small $^4$  Departamento de Matem\'atica, FCE}\\[-1mm] {\it \small Universidad Austral}\\[-1mm]
		{\it \small Paraguay 1950, S2000FZF Rosario, Argentina}\\
	{\it\small $^5$  CONICET, Argentina}}}


\date{}
\maketitle
\thispagestyle{empty}

\vskip 0mm

\noindent {\small{\bf Abstract.}
We consider an elliptic variational inequality with unilateral constraints in a Hilbert space $X$ which,
under appropriate assumptions on the data, has a unique solution $u$.
We formulate a convergence criterion to the solution $u$, i.e., we provide necessary and sufficient conditions on a sequence $\{u_n\}\subset X$ which guarantee  the convergence $u_n\to u$ in the space $X$. Then, we illustrate the use of this criterion to recover well-known convergence results and well-posedness results in the sense of Tykhonov and Levitin-Polyak.  We also provide two applications of our results, in the study of a heat transfer problem and an elastic frictionless contact problem, respectively.}

\vskip 3mm
\noindent
{\bf Keywords :}	Elliptic variational inequality,  convergence criterion, convergence results, well-posedness, contact, heat transfer,
unilateral constraint.

\vskip 3mm
\noindent {\bf 2010 MSC:} \ 47J20, 49J40, 40A05, 74M15, 74M10, 35J20.


\section{Introduction}\label{s1}
\setcounter{equation}0


\medskip\noindent
A large number of mathematical models in Physics, Mechanics and Engineering Science are expressed in terms of
strongly  nonlinear boundary value problems for partial  differential equations which, in a weak formulation, lead to variational inequalities.
The theory of variational inequalities was developed  based on
arguments of monotonicity and convexity, including properties of the subdifferential
of a convex function.  Because
of their importance in partial differential equations theory and engineering applications, a considerable effort has been put into the analysis, the control
and the numerical simulations of variational inequalities.
Basic references in the field are
\cite{BC, B, G, Kind-St,Li}, for instance.
Applications of variational inequalities in Mechanics can be found in
the books \cite{C,DL,EJK,HS,HHNL,KO, P}.

In this paper we study the convergence of an arbitrary sequence to the solution of an elliptic variational inequality.   Our results below could be extended to more general inequalities in reflexive Banach spaces. Nevertheless, for simplicity, we restrict ourselves to the following functional framework:   $X$ is a real Hilbert space endowed with the inner product $(\cdot,\cdot)_X$  and the associated norm  $\|\cdot\|_X$,
$K\subset X$, $A:X\to X$, $j:X\to\R$ and $f\in X$. Then, the inequality problem we consider in this paper is as follows.

\medskip\medskip\noindent{\bf Problem}  ${\cal P}$. {\it Find $u$ such that}
\begin{equation}\label{1}u\in K,\qquad(Au,v-u)_X+j(v)-j(u) \ge(f,v-u)_X \qquad\forall\,v\in K.
\end{equation}

The unique  solvability of Problem $\cP$
follows from well-known results obtained in the literature, under various assumptions on the data. Here, we shall use the existence and uniqueness results that we recall in the next section, Theorem \ref{t1}. We also present some convergence results to the solution $u$ of inequality (\ref{1}). Note that  a large number of convergence results for inequality \eqref{1} have been obtained in literature.   The continuous dependence of the solution with respect to the data, the convergence of the solution to  penalty problems when the penalty parameter converges to zero, the convergence of the solutions of discrete numerical schemes, the convergence of the solution of various perturbed problems when some parameters converge are several examples, among others.  Note also that the concept of well-posedness (in the sense of Tyknonov or Levitin-Polyak) for inequality (\ref{1})  is also based on the convergence to the solution  $u$ of the so-called approximating and generalized approximating sequences, respectively.

All these examples, together with various relevant applications in Optimal Control Theory, Physics and Mechanics, lead to the following question: is it possible to describe the convergence of a sequence of  $\{u_n\}\subset X$ to the solution $u$ of the variational inequality (\ref{1})? In other words, the question is to provide necessary and sufficient conditions for the convergence $u_n\to u$ in $X$, i.e., to  provide a {\it convergence criterion}. The first aim of this paper is to provide an answer to this question. Here, we state and prove such a criterion of convergence expressed in terms of metric properties. The second aim is to illustrate how this  criterion could be used in various examples and applications, in order to deduce some convergence results.

A short description of the rest of the manuscript is as follows.
First, in Section~\ref{s2} we present preliminary results concerning the unique solvability of Problem $\cP$, together with some convergence results. In Section~\ref{s3} we state and prove our main result,  Theorem~\ref{tm}, which represents a criterion of convergence to the solution $u$ of inequality~(\ref{1}).
Section~\ref{s4} is devoted to some theoretical applications  of Theorem~\ref{tm}.
Here we state and prove two convergence results,  introduce a new well-posedness concept and show that it extends the classical
Tykhonov and Levitin-Polyak well-posedness concepts for variational inequalities of the form~(\ref{1}).
Finally, in Sections~\ref{s5}  and~\ref{s6}  we present applications  of our theoretical results in the study of two specific boundary value problems, which model  a static frictionless contact process for elastic materials and
a stationary {heat transfer} problem,  respectively. For these problems we state and prove convergence results and provide their physical and mechanical interpretations.

\section{Preliminaries}\label{s2}
\setcounter{equation}0

Everywhere in this paper, unless it is specified otherwise, we use the functional framework described in Introduction. Notation $0_X$ and $I_X$ will represent the zero element and the identity operator of $X$, respectively. All the limits, upper and lower limits below are considered as $n\to\infty$, even if we do not mention it explicitly. The symbols ``$\rightharpoonup$"  and ``$\to$"
denote the weak and the strong convergence in various spaces which will be specified, except in the case when these convergence take place in $\real$.  For a sequence $\{\ve_n\}\subset\R_+$ that converges to zero we use the short hand notation $0\le\ve_n\to 0$. Finally,
we denote by $d(u,K)$ the distance between an element $u\in X$ to the set $K$, that is
\begin{equation}\label{d}
d(u,K)=\inf_{v\in K} \|u-v\|_X.
\end{equation}

For the convenience of the reader we also recall the following definitions { which can be found in many books and surveys, including \cite{DMP,Z}, for instance.
}

\begin{Definition}\label{Mosco}
	Let $\{K_n\}$ be a sequence of nonempty subsets of \ $X$ and let $K$ be a  nonempty subset of $X$.
	We say that the sequence $\{K_n\}$ converges to $K$ {\it in the sense of Mosco} (\cite{Mosco}) and we write $K_n\stackrel{M}{\longrightarrow} K$, if the following conditions hold:
	
	\smallskip
	{\rm 	(a)} for each $u\in K$, there exists a sequence $\{u_n\}$ such that $u_n\in K_n$ for each $n\in \mathbb{N}$ and $u_n\rightarrow u$ in $X$.
	
	\smallskip	
	{\rm  (b)} for each sequence $\{u_n\}$ such that $u_n\in K_n$ for each $n\in \mathbb{N}$ and  $u_n\rightharpoonup u$  in $X$, we have $u\in K$.
	
\end{Definition}

\begin{Definition}\label{AopH}
	{ An operator $A \colon X \to X$ is called:
		
		\smallskip
		${\rm(a)}$\,
		{\em monotone} if $(Au - A v, u-v)_X \ge 0$ \ $\forall\, u, v\in X$.

		\smallskip
		${\rm(b)}$
		{\em strongly monotone} if there exists $m>0$ such that
		\begin{equation}\label{A1}
			(Au-Av,u-v)_X\ge m\|u-v\|_X^2 \ \quad\forall\,
			u,v\in X.
		\end{equation}

		\smallskip
		${\rm(c)}$
		{\em pseudomonotone},
		if it is bounded and  the convergence $u_n \rightharpoonup u$ in $X$ together with inequality\
		$\displaystyle \limsup\, (A u_n, u_n -u)_X \le 0$
		imply that}
	\[\displaystyle
	\liminf\, (A u_n, u_n - v)_X\ge (A u, u - v)_X\quad\ \forall\,v \in X.\]
	
	\smallskip
	${\rm(d)}$
	{\em hemicontinuous}
	if for all $u$, $v$, $w \in X$, the function
	$\lambda \mapsto (A (u + \lambda v), w )_X$
	is continuous on $[0, 1]$.
	
	\smallskip
	${\rm(e)}$
	{\em demicontinuous}
	if $u_n \to u$ in $X$ implies $A u_n\rightharpoonup Au$ in $X$.

	\smallskip
	${\rm(d)}$
	{\em Lipschitz continuous}
	if there exists $M>0$ such that
	\begin{equation}
	\label{A2}\|Au-Av\|_X\le M\,\|u-v\|_X
	\quad\forall\,
	u,v\in X.
	\end{equation}

\end{Definition}

{ Next, we follow \cite[p.267]{PS} and introduce the following notion of 
a penalty operator.}

\begin{Definition}\label{penalty}
	An operator $G:X\to X$ is said to be a penalty operator of $K$ if it is  bounded, demicontinuous, monotone, and $Gu=0_X$ if and only if $u\in K$.
\end{Definition}

Note that if $K$ { is a nonempty, closed, convex subset of $X$ and $P_K$ denotes the projection operator on $K$, then it is easy to see that the operator $G=I_X-P_K:X\to X$ is monotone and Lipschitz continuous. Therefore, using Definitions \ref{AopH} and \ref{penalty} it follows that $G$ is a penalty operator of $K$.} Moreover,  the proposition below, proved in \cite{Z}, shows that any penalty operator is pseudomonotone.

\begin{Proposition}\label{1pE1} Let $G:X\to X$ be a bounded, hemicontinuous and monotone operator. Then $G$ is pseudomonotone.
\end{Proposition}

In addition, the following result, stated and proved in \cite{Z}, concerns the sum of two pseudomonotone operators.

\begin{Proposition}\label{2pE1} Let $A,\ B :X\to X$ be pseudomonotone operators. Then the sum $A+B:X\to X$ is a pseudomonotone operator.
\end{Proposition}

\medskip
In the study of Problem $\cP$ we consider the following assumptions.
\begin{eqnarray}
&&\label{K}
K \ \mbox{\rm is a nonempty, closed, convex subset of} \ X.
\\[2.5mm]
&&\label{Ap}
\left\{ \begin{array}{l}
A:X\to X\ {\rm is\ a\ pseudomonotone\ operator\ and}\\[0mm]
{\rm there\ exists}\ m>0\ \ {\rm such\ that\ \eqref{A1}\ holds.}
\end{array}\right.\\[2.5mm]
&&\label{j}j:X\to \mathbb{R} \mbox{ is
	convex  and lower semicontinuous. }\\[2mm]
&&\label{f}f\in X.
\end{eqnarray}

\medskip
We now recall  the following  well-known existence and uniqueness result.

\begin{Theorem}\label{t1}
	Assume $(\ref{K})$--$(\ref{f})$.
	Then, the variational inequality $(\ref{1})$  has a unique solution $u$.
\end{Theorem}

Theorem \ref{t1} represents a particular case of Theorem 84 in \cite{SMBOOK}. We now complete it with some convergence results of the form
\begin{equation}\label{c0}
u_n\to u\quad{\rm in}\quad X.
\end{equation}
Here and everywhere in Sections \ref{s2} and \ref{s3} we keep notation $u$ for the unique solution of Problem $\cP$ obtained in Theorem \ref{t1}, even if we do not mention explicitly. Moreover, $\{u_n\}$ represents a sequence of elements of $X$ which will be specified.

Consider the sequences $\{K_n\}$, $\{\lambda_n\}$, $\{f_n\}$ such that the following conditions hold, for each $n\in \mathbb{N}$.
\begin{eqnarray}
&&\label{Kn}
K_n \ \mbox{\rm is a nonempty, closed, convex subset of} \ X.\\ [2mm]
&&\label{la}\lambda_n>0.\\ [2mm]
&&\label{fn}f_n\in X.
\end{eqnarray}
Assume also that
\begin{equation}\label{G}G:X\to X \ \ \mbox{is a  penalty operator of}\ K.
\end{equation}
Then, using Propositions  \ref{1pE1} and  \ref{2pE1} it follows that the operator
$A+\frac{1}{\lambda_n} G:X\to X$ satisfies condition (\ref{Ap}), for each $n\in\mathbb{N}$.
Therefore, under the previous assumptions,
Theorem \ref{t1} guarantees the unique solvability of the following three problems.

\medskip\medskip\noindent{\bf Problem}  ${\cal P}_n^1$. {\it Find $u_n$ such that}
\begin{equation}
\label{1K}u_n\in K_n,\quad(Au_n,v-u_n)_X+j(v)-j(u_n) \ge(f,v-u_n)_X \quad\forall\,v\in K_n.
\end{equation}

\medskip\medskip\noindent{\bf Problem}  ${\cal P}_n^2$. {\it Find $u_n$ such that}
\begin{eqnarray}
&&\label{1G}u_n\in X,\quad(Au_n,v-u_n)_X+\frac{1}{\lambda_n}(Gu_n,v-u_n)_X+j(v)-j(u_n)\nonumber\\ [0mm]
	&&\qquad\qquad\qquad\ge(f,v-u_n)_X \quad\forall\,v\in X.
\end{eqnarray}	
	
\medskip\medskip\noindent{\bf Problem}  ${\cal P}_n^3$. {\it Find $u_n$ such that}	
\begin{equation}\label{1f}u_n\in K,\quad(Au_n,v-u_n)_X+j(v)-j(u_n) \ge(f_n,v-u_n)_X \quad\forall\,v\in K.
\end{equation}

\medskip
Moreover, we have the following result.

\medskip
\begin{Theorem}\label{t2}
Assume	 $(\ref{K})$--$(\ref{f})$.
Then, the convergence $(\ref{c0})$ holds in each of the following three cases.
	
\medskip	
{\rm a)} Condition $(\ref{Kn})$ holds,  $K_n\stackrel{M}{\longrightarrow} K$, and $u_n$ denotes the solution to Problem ${\cal P}_n^1$.

\medskip	
{\rm b)} Conditions  $(\ref{la})$ and  $(\ref{G})$ hold, $\lambda_n\to 0$,
and $u_n$ denotes the solution to Problem ${\cal P}_n^2$.

\medskip	
{\rm c)} Condition $(\ref{fn})$ holds,  $f_n\to f$ in $X$,
and $u_n$ denotes the solution to Problem ${\cal P}_n^3$.
\end{Theorem}

Theorem \ref{t2} represents a version of some convergence results obtained in \cite{SofMat, ST2}, for instance and, for this reason, we skip  its proof.
{ Nevertheless, we mention that a proof of the convergence results in Theorem \ref{t2} b), c) will be provided in Section \ref{s4},
under additional assumptions on the operator $A$ and the function $j$. This proof is based on the convergence criterion we state and prove in Section \ref{s3}. Moreover, for the convenience of the reader, we present below a sketch of the proof of the point a) of this theorem, based
on standard arguments of compactness, monotonicity, pseudomonotonicity and lower semicontinuity.

\medskip\noindent
{\it Proof of Theorem  $\ref{t2}\,a)$}. The proof is structured in   three steps, as follows.

\smallskip
{\it Step i)} We use inequality \eqref{1} and the strong monotonicity of the operator $A$ in order to prove that the sequence $\{u_n\}$ is bounded in $X$. Then, using the reflexivity of  the space $X$, we deduce that this sequence contains a subsequence, again denoted by $\{u_n\}$,  such that $u_n\rightharpoonup\widetilde{u}$ with some
$\widetilde{u}\in X$.

\smallskip
{\it Step ii)} We use the pseudomonotonicity of $A$ and the properties of the function $j$ to deduce that the element $\widetilde{u}$ satisfies inequality \eqref{1}. Therefore, by the uniqueness of the solution, we deduce that $\widetilde{u}=u$. Moreover, a careful analysis reveals that any weakly convergent subsequence of the sequence   $\{u_n\}$ converges weakly to $u$, in $X$.
We then use a standard argument to deduce that the whole sequence $\{u_n\}$ converges weakly in $X$ to~$u$.

\smallskip
{\it Step iii)} Finally, we use the strong monotonicity of the operator $A$, and the weak convergence  $u_n\rightharpoonup u$ in $X$ to deduce that the strong convergence holds, $(\ref{c0})$, which concludes the proof.
\hfill$\Box$}

\medskip
We end this section with the remark that Theorem \ref{t2} provides relevant examples of sequences which converge to the solution $u$ of the variational inequality (\ref{1}). However, some elementary examples show that, besides the sequences introduced in  parts a), b) and c) of this theorem, there  exists other sequences $\{u_n\}$ which converge to $u$. A criterion which could identify all such sequences is provided in the next section.

\section{A convergence criterion}\label{s3}
\setcounter{equation}0

We now consider the following additional condition on the operator $A$ and the function~$j$:
\begin{eqnarray}
&&\hspace{-11mm}\label{A}
\left\{ \begin{array}{l}
A:X\to X\ {\rm is\ a\ strongly\ monotone\ Lipschitz\ continuous\ operator,}\\[0mm]
{\rm i.e.,\ there\ exist}\ m>0\ {\rm and}\ M>0\ {\rm such\ that\ \eqref{A1}\ and\ \eqref{A2}\ hold.}
\end{array}\right.\\[2.5mm]
&&\hspace{-11mm}\label{jj}
\left\{\begin{array}{ll} j:X\to \mathbb{R}\ \mbox{ is
convex  and  for each $D>0$ there exists $L_D>0$ such that}\\ [2mm]
|j(u)-j(v)|\le L_D \|u-v\|_X\quad\forall\, u,\ v\in X\ {\rm with}\ \|u\|_X, \ \|v\|_X\le D.
\end{array}\right.	
\end{eqnarray}
Note that conditions \eqref{A} and \eqref{jj} imply \eqref{Ap} and \eqref{j}, respectively. Therefore, Theorems \ref{t1} and \ref{t2} still hold under assumptions $(\ref{K})$, $(\ref{f})$, $(\ref{A})$ and  $(\ref{jj})$. Moreover,   condition  \eqref{jj} shows that the function $j$ is Lipschitz continuous on each bounded sets of $X$. This implies that, in particular,  $j$ is locally Lipschitz. In addition, note that any continuous seminorm on the space $X$ satisfies condition \eqref{jj}.

\medskip
Our main result in this section is the following.

\begin{Theorem}\label{tm}
	Assume $(\ref{K})$, $(\ref{f})$, $(\ref{A})$, $(\ref{jj})$ { and denote by $u$ the solution of the variational inequality $(\ref{1})$ provided by Theorem $\ref{t1}$. Consider also an arbitrary sequence
	$\{u_n\}\subset X$, together with the statements  \eqref{c1} and \eqref{c2} below.}
\begin{eqnarray}
	&&\label{c1}
	\quad\ \ u_n\to u\qquad{\rm in}\ X.
	\\ [4mm]
	&&\label{c2}
		\left\{\begin{array}{ll}  {\rm (a)}\ d(u_n,K)\to 0\ ;\\[4mm]
		 {\rm (b)}\ 	\mbox{\rm there exists $0\le\ve_n\to 0$ such that} \\[3mm]	
\quad\quad(Au_n,v-u_n)_X+j(v)-j(u_n)+\ve_n(1+\|v-u_n\|_X)\\ [1mm]
\qquad\qquad\qquad\quad \ge(f,v-u_n)_X  \quad\forall\,v\in K,\ n\in\mathbb{N}.
\end{array}\right.	
\end{eqnarray}
{ Then, these statements are equivalent, i.e.,  \eqref{c1} holds  if and only if \eqref{c2} holds.}	
\end{Theorem}

\medskip
The proof of Theorem \ref{tm} is based on the following result.

\begin{Lemma}\label{lm}
	Assume $(\ref{K})$, $(\ref{f})$, $(\ref{A})$, $(\ref{jj})$. Then  any sequence $\{u_n\}\subset X$  which satisfies condition $(\ref{c2})\,({\rm b)}$ is bounded.
\end{Lemma}

\medskip\noindent
{\it Proof.} Let $n\in\mathbb{N}$. We test in
(\ref{c2}) with a fixed element of $K$, say $v=u$. We have
\begin{equation*}
(Au_n,u-u_n)_X+j(u)-j(u_n)+\ve_n(1+\|u-u_n\|_X)\ge(f,u-u_n)_X  	
\end{equation*}	
which implies that
\begin{eqnarray*}
&&(Au_n-Au,u_n-u)_X+j(u_n)\\ [2mm]
&&\qquad\le (Au,u-u_n)_X+ j(u)+\ve_n(1+\|u-u_n\|_X)+(f,u_n-u)_X  	\nonumber
\end{eqnarray*}
and, moreover,
\begin{eqnarray}
	&&\label{rd1}m\,\|u-u_n\|^2_X+j(u_n)\\ [2mm]
	&&\qquad\le \|Au\|_X\|u-u_n\|_X+ j(u)+\ve_n(1+\|u-u_n\|_X)+\|f\|_X\|u_n-u\|_X.  	\nonumber
\end{eqnarray}
On the other hand,  assumption (\ref{jj}) {and a standard result on convex lower semicontinuous functions (see \cite[p.208]{KZ}, for instance}) implies that $j$ is bounded from below by an affine continuous function. Hence, there exist $\alpha\in X$, $\beta\in \mathbb{R}$ such that
\begin{equation}\label{zz}
j(v)\ge (\alpha,v)_X+\beta\qquad\forall \,v\in X.	
\end{equation}
This implies that
\begin{eqnarray*}
&&j(u_n)\ge (\alpha,u_n)_X+\beta=(\alpha,u_n-u)_X+(\alpha,u)_X+
\beta\\ [2mm]
&&\qquad
\ge-\|\alpha\|_X\|u_n-u\|_X-\|\alpha\|_X\|u\|_X-|\beta|
\end{eqnarray*}
and, substituting this inequality in (\ref{rd1})
yields
\begin{eqnarray*}
&&m\,\|u-u_n\|^2_X\le (\|\alpha\|_X+\|Au\|_X+\|f\|_X+\ve_n)\|u-u_n\|_X\\ [2mm]
&&\qquad\qquad+j(u)+\|\alpha\|_X\|u\|_X+|\beta|+\ve_n.  	
\end{eqnarray*}
{Moreover, since $j(u)\le |j(u)|$, we deduce that
\begin{eqnarray*}
	&&m\,\|u-u_n\|^2_X\le (\|\alpha\|_X+\|Au\|_X+\|f\|_X+\ve_n)\|u-u_n\|_X\\ [2mm]
	&&\qquad\qquad+|j(u)|+\|\alpha\|_X\|u\|_X+|\beta|+\ve_n.  	
\end{eqnarray*}}
Next, we use  the implication
\begin{equation}\label{in}
x^2\le ax+b\quad \Longrightarrow\quad x\le a+\sqrt{b}\qquad\forall\, x,\, a,\ b\ge 0
\end{equation}
and the convergence $\ve_n\to 0$ to see that the sequence $\{\|u-u_n\|_X\}$ is bounded in $\mathbb{R}$. This implies that the sequence $\{u_n\}$ is bounded in $X$, which concludes the proof.
\hfill$\Box$

\medskip
We now proceed with the proof of Theorem \ref{tm}.

\medskip\noindent
{\it Proof.}  Assume that \eqref{c1} holds. Then, since $u\in K$ it follows that $d(u_n,K)\le \|u_n-u\|_X$ for each $n\in\mathbb{N}$ which implies that  (\ref{c2})(a) holds. To prove (\ref{c2})(b) we fix  $n\in\mathbb{N}$ and $v\in K$. We write
\begin{eqnarray*}&&(Au_n,v-u_n)_X+j(v)-j(u_n)-(f,v-u_n)_X\\ [2mm]
&&\quad=(Au_n-Au,v-u_n)_X+(Au,v-u)_X+(Au,u-u_n)_X\\ [2mm]
&&\qquad+j(v)-j(u)+j(u)-j(u_n)-(f,v-u)_X+(f,u_n-u)_X
\end{eqnarray*}
and, using \eqref{1} we deduce that
\begin{eqnarray}&&\label{r23}(Au_n,v-u_n)_X+j(v)-j(u_n)-(f,v-u_n)_X\\ [2mm]
	&&\quad\ge(Au_n-Au,v-u_n)_X+(Au,u-u_n)_X +j(u)-j(u_n)+(f,u_n-u)_X.\nonumber
\end{eqnarray}
We now use \eqref{r23} and  inequalities
\begin{eqnarray*}
&&(Au_n-Au,v-u_n)_X\ge -\|Au_n-Au\|_X\|v-u_n\|_X\ge-M\|u-u_n\|_X\|v-u_n\|_X,\\ [2mm]
&&(Au,u-u_n)_X\ge-\|Au\|_X\|u-u_n\|_X,\\[2mm]
&&(f,u_n-u)_X\ge -\|f\|_X\|u-u_n\|_X
\end{eqnarray*}
to find that	
\begin{eqnarray*}&&(Au_n,v-u_n)_X+j(v)-j(u_n)-(f,v-u_n)_X+j(u_n)-j(u)\\ [2mm]
&&\quad+M\|u-u_n\|_X\|v-u_n\|_X+\|Au\|_X\|u-u_n\|_X+
\|f\|_X\|u-u_n\|_X\ge 0.
\end{eqnarray*}
Therefore, with notation
\begin{equation}\label{ep}\ve_n={\rm max}\,\{M\|u-u_n\|_X,(\|Au\|_X+
	\|f\|_X)\|u-u_n\|_X+j(u_n)-j(u)\}
\end{equation}
we see that
\begin{equation}\label{e1n}
(Au_n,v-u_n)_X+j(v)-j(u_n)+\ve_n(1+\|v-u_n\|_X)\ge(f,v-u_n)_X.
\end{equation}
On the other hand, notation \eqref{ep}, assumption \eqref{c1} and
the continuity of the function $j$, guaranteed by hypothesis \eqref{jj}, show that
\begin{equation}\label{e2n}
\ve_n\to 0.
\end{equation}
We now combine (\ref{e1n}) and (\ref{e2n}) to see that condition (\ref{c2})(b) is satisfied.

Conversely, assume now that (\ref{c2}) holds.
Then, (\ref{c2})(a) and  definition \eqref{d} of the distance function show that for each $n\in\mathbb{N}$ there exist two elements $v_n$ and $w_n$ such that
\begin{equation}\label{e3}
u_n=v_n+w_n,\quad v_n\in K,\quad w_n\in X,\quad \|w_n\|_X\to 0.
\end{equation}
We fix $n\in\mathbb{N}$ and use (\ref{c2})(b)
with $v=u\in K$ to see that
\begin{equation}\label{e4}
	(Au_n,u-u_n)_X+j(u)-j(u_n)+\ve_n(1+\|u-u_n\|_X)\ge(f,u-u_n)_X.
\end{equation}
On the other hand, we use the regularity $v_n\in K$ in (\ref{e3}) and test with $v=v_n$ in \eqref{1} to find that
\begin{equation}\label{e5}
	(Au,v_n-u)_X+j(v_n)-j(u)\ge(f,v_n-u)_X.
\end{equation}
We now add inequalities (\ref{e4}), (\ref{e5})  to obtain that
\begin{eqnarray}
&&\label{Z}(Au_n,u-u_n)_X+(Au,v_n-u)_X
	+j(v_n)-j(u_n)\\ [2mm]
&&\quad+\ve_n(1+\|u-u_n\|_X)\ge(f,v_n-u_n)_X.\nonumber
\end{eqnarray}
Next, we use equality $u_n=v_n+w_n$ to see that
\begin{eqnarray*}
	&&(Au_n,u-u_n)_X+(Au,v_n-u)_X=
	(Au,v_n-u)_X-(Av_n,v_n-u)_X\\ [2mm]
	&&\quad+(Av_n,v_n-u)_X+(Au_n,u-v_n)_X-(Au_n,w_n)_X
	\\ [2mm]
	&&\qquad=(Au-Av_n,v_n-u)_X+(Au_n-Av_n,u-v_n)_X-(Au_n,w_n)_X
\end{eqnarray*}	
and, therefore, \eqref{Z} implies that
	\begin{eqnarray*}
	&&(Au-Av_n,v_n-u)_X+(Au_n-Av_n,u-v_n)_X-(Au_n,w_n)
	+j(v_n)-j(u_n)\\ [2mm]
	&&\quad+\ve_n(1+\|u-v_n-w_n\|_X)+ (f,w_n)_X\ge 0.
\end{eqnarray*}
Using now  assumption \eqref{A} and equality $u_n=v_n+w_n$ we deduce that
\begin{eqnarray}
&&\label{e6}m\|u-v_n\|^2_X\le M\|w_n\|_X\|u-v_n\|_X+ \|Au_n\|_X
\|w_n\|_X\\ [2mm]
&&\qquad+j(v_n)-j(u_n)+\ve_n+\ve_n\|u-v_n\|_X
+\ve_n\|w_n\|_X+\|f\|_X\|w_n\|_X.\nonumber
\end{eqnarray}

On the other hand, Lemma  \ref{lm} guarantees that the sequence $\{u_n\}$ is bounded in $X$. Therefore, \eqref{e3} implies that there exists $D>0$ such that
\begin{equation}\label{e10}
\|u_n\|_X\le D,\qquad \|v_n\|_X\le D\qquad\forall\, n\in\mathbb{N}.
\end{equation}
Moreover, using condition (\ref{A}) we may assume that
\begin{equation}\label{e9}\|Au_n\|_X\le D\qquad\forall\, n\in\mathbb{N}.
\end{equation}
In addition, using  (\ref{e10}),  assumption (\ref{jj}) and equality $u_n=v_n+w_n$ in \eqref{e3} we find that
\begin{equation}\label{e7}
j(v_n)-j(u_n)\le L_D\|w_n\|_X.
\end{equation}

We now combine the bounds \eqref{e6},   \eqref{e9} and \eqref{e7} to deduce that
\begin{eqnarray}
	&&\label{e8}m\|u-v_n\|^2_X\le (M\|w_n\|_X+\ve_n)\|u-v_n\|_X\\ [2mm]
	&&\qquad+ (D
	+L_D+\ve_n+\|f\|_X)\|w_n\|_X+\ve_n.\nonumber
\end{eqnarray}
Next, we use (\ref{e8}),  inequality (\ref{in})
and the convergences $\|w_n\|_X\to 0$, $\ve_n\to 0$ to find that
$\|u-v_n\|_X\to 0$. This implies that $v_n\to u$ in $X$ and, using \eqref{e3} we deduce that
\eqref{c1} holds, which concludes the proof.\hfill$\Box$

\medskip
Theorem \ref{tm} shows that, under assumptions
$(\ref{K})$--$(\ref{f})$, $(\ref{A})$, $(\ref{jj})$, conditions $(\ref{c2})({\rm a})$ and  $(\ref{c2})({\rm b})$ represent necessary and sufficient conditions for the convergence $(\ref{c1})$. The two examples below show that, in general, we cannot skip one of these conditions.

\begin{Example} Take $A=I_X$ and $j\equiv 0$. Then the solution of inequality $(\ref{1})$ is $u=P_Kf$ where, recall, $P_K$ represents the projector operator on $K$. Assume now that $f\notin K$ and take $u_n=f$ for each $n\in\mathbb{N}$. It follows from here that  $(\ref{c2})({\rm b})$ is satisfied with $\ve_n=0$.
Nevertheless,  $(\ref{c2})({\rm a})$ does not hold since
\[d(u_n,K)=\|u_n-P_Ku_n\|_X=\|f-P_Kf\|_X>0.\]
Moreover, the convergence   $(\ref{c1})$ is not valid. We conclude from here that
condition $(\ref{c2})({\rm a})$ cannot be skipped, i.e., condition  $(\ref{c2})({\rm b})$ is not a sufficient condition to guarantee the convergence \eqref{c1}.
\end{Example}

\begin{Example}
Let $X=\real$, $K=[0,1]$, $A=I_X$, $j\equiv 0$ and  $f=\frac{1}{2}$. Then, it follows that $u=f=\frac{1}{2}$ and the sequence $\{u_n\} $ with $u_n=0$ satisfies $(\ref{c2})({\rm a})$  but does not satisfy $(\ref{c1})$. We conclude that
condition  $(\ref{c2})({\rm a})$ is not a sufficient condition to guarantee
the convergence $(\ref{c1})$.
\end{Example}

We end this section with the remark that Theorem \ref{tm} was obtained under the assumptions (\ref{A}), (\ref{jj}) which, however, are not necessary neither in the statement of Theorem \ref{t1} nor in the statement of Theorem \ref{t2}. Removing or relaxing these conditions is an interesting problem which clearly deserves to be investigated into future.

\section{Convergence and well-posedness results}\label{s4}
\setcounter{equation}0

In this section we show how Theorem \ref{tm} can be used to deduce some theoretical convergence results.
We also prove that the well-posedness of inequality  (\ref{1}), both in the Tykhonov and Levitin-Polyak sense, can be deduced as a consequence of this theorem. Finally
we provide an interpretation of Theorem \ref{tm}
in the context of the $\cT$-well-posedness concept introduced in \cite{S,SX16} and used in various papers, including \cite{ST1,ST3}.

\medskip\noindent
{\bf Convergence results.} A first consequence of Theorem \ref{tm} is the  following continuous dependence result.

\begin{Corollary}\label{cor1}
	Assume	 $(\ref{K})$, $(\ref{f})$, $(\ref{fn})$, $(\ref{A})$  and $(\ref{jj})$. Also, assume that
	$f_n\to f$ in $X$ and denote by $u_n$ the solution  of Problem ${\cal P}_n^3$.
	Then, the convergence $(\ref{c0})$ holds.
\end{Corollary}

\medskip\noindent{\it Proof.}  Let $n\in\mathbb{N}$ and $v\in K$. We use \eqref{1f} and write
\begin{equation*}u_n\in K,\quad(Au_n,v-u_n)_X+j(v)-j(u_n)+(f-f_n,v-u_n)_X \ge(f,v-u_n)_X,
\end{equation*}
which implies that
\begin{equation*}u_n\in K,\quad(Au_n,v-u_n)_X+j(v)-j(u_n)+\|f-f_n\|_X\|v-u_n\|_X \ge(f,v-u_n)_X.
\end{equation*}
It follows from this inequality  that conditions (\ref{c2})(a) and  (\ref{c2})(b) are satisfied with $\ve_n=\|f-f_n\|_X\to 0$. We now use Theorem \ref{tm} to conclude the proof.
\hfill$\Box$

\medskip
A second consequence of Theorem \ref{tm} concerns the penalty problem $\cP_n^2$ and is as follows.

\begin{Corollary}\label{cor2}
	Assume	 $(\ref{K})$, $(\ref{f})$,
	$(\ref{la})$, $(\ref{G})$, $(\ref{A})$ and $(\ref{jj})$. Also, assume that
	$\lambda_n\to 0$ and denote by $u_n$ the solution  of Problem ${\cal P}_n^2$.
	Then, the convergence $(\ref{c0})$ holds.
\end{Corollary}

\medskip\noindent{\it Proof.} The proof is structured in several steps, as follows.

\medskip
\noindent {\it Step {i)}  We prove that the sequence $\{u_n\}$ satisfies condition $(\ref{c2}){\rm (b)}$.}
Let $n\in\mathbb{N}$ and $v\in K$. Then,  (\ref{G}) and Definition $\ref{penalty}$ imply that $Gv=0_X$, $(Gv-Gu_n,v-u_n)_X\ge 0$
and, therefore,
\begin{equation}\label{zzz}
(Gu_n,v-u_n)_X\le 0.
\end{equation}
We now use inequalities (\ref{1G}) and (\ref{zzz}) to see that
\begin{equation}\label{m4}
(Au_n,v-u_n)_X+j(v)-j(u_n)\ge(f,v-u_n)_X
\end{equation}
which shows that ${\rm (\ref{c2})(b)}$ holds with $\ve_n=0$.

\medskip
\noindent {\it Step {ii)}  We prove that any weakly convergent subsequence of the sequence $\{u_n\}$  satisfies condition ${\rm (\ref{c2})(a)}$.} Indeed,
consider a weakly convergent subsequence of the sequence $\{u_n\}$, again denoted by $\{u_n\}$.
Then, there exists an element $\widetilde{u}\in X$ such that
\begin{equation*}
u_n\rightharpoonup \widetilde{u}\quad{\rm in}\quad X.
\end{equation*}	
We shall prove that $\widetilde{u}\in K$ and $u_n\to \widetilde{u}$ in $X$. To this end, we fix $n\in\mathbb{N}$ and $v\in X$. We use (\ref{1G}) and (\ref{zz}) to write	
\begin{eqnarray*}
		&&\frac{1}{\lambda_n}( G{u}_n,{u}_n-v)_X\leq ( A{u}_n,v-{u}_n)_X+j(v)-j({u}_n)+(f,{u}_n-v)_X\\ [2mm]
		&&\quad \leq \|Au_n\|_{X}\|v-u_n\|_X+j(v)+\|\alpha\|_{X}\|u_n\|_{X}+|\beta|+
		\|f\|_{X}\|v-u_n\|_X.
\end{eqnarray*}

On the other hand, Step i) and Lemma \ref{lm} imply that the sequence $\{u_n\}$ is bounded. Therefore,  from the previous inequality
we deduce that there exists a constant $C(v)>0$ which does not depend on $n$ such that
	\begin{equation*}
	(G{u}_n,{u}_n-v)_X\leq\lambda_nC(v).
	\end{equation*}
	Passing to the upper limit in the above inequality
	and using assumption $\lambda_n\to 0$ we find that
	\begin{equation}\label{Ee4.14}
		\lim\sup\,(G{u}_n,{u}_n-v)_X\leq 0.
	\end{equation}
	Taking now $v=\widetilde{u}$ in (\ref{Ee4.14}) we deduce that
	\[\limsup\,( G{u}_n,{u}_n-\widetilde{u})_X\leq 0.\]
	Then, using the convergence ${u}_n \rightharpoonup {\widetilde{u}}$ in $X$ and the pseudomonotonicity of $G$, guaranteed by assumption (\ref{G}) and Proposition \ref{1pE1}, we deduce that
	\begin{eqnarray}\label{Ee4.16}
		( G\widetilde{u},\widetilde{u}-v)_X\leq \liminf\,( G{u}_n,{u}_n-v)_X.
	\end{eqnarray}
	We now combine inequalities (\ref{Ee4.14}) and (\ref{Ee4.16}) to see that
	
	\begin{equation*}
		( G\widetilde{u},\widetilde{u}-v)_X\leq 0.
	\end{equation*}
	Recall that this inequality  is valid for any $v\in X$. Therefore, we deduce that
	$G\widetilde{u}=0_X$ and, using assumption
	(\ref{G}) we find that
	\begin{equation}\label{Ereg}
		\widetilde{u}\in K.
	\end{equation}
	
	Let $n\in\mathbb{N}$. Then, using (\ref{m4}) with $v=\widetilde{u}$ we find that
	\begin{equation*}
		( A{u}_n,{u}_n-\widetilde{u})_X
		\le j(\widetilde{u})-j(u_n)+( f,{u}_n-\widetilde{u})_X
	\end{equation*}
or, equivalently,
\begin{equation*}
	(A{u}_n-A\widetilde{u},{u}_n-\widetilde{u})_X
	\le (A\widetilde{u},\widetilde{u}-{u}_n)_X+j(\widetilde{u})-j(u_n)+( f,{u}_n-\widetilde{u})_X.
\end{equation*}
Finally, we use the strong monotonicity of the operator $A$, (\ref{A1}), to see that
\[
m\|u_n-\widetilde{u}\|_X^2\le
(A\widetilde{u},\widetilde{u}-{u}_n)_X+j(\widetilde{u})-j(u_n)+( f,{u}_n-\widetilde{u})_X.\]

	We now pass  to the upper limit in this inequality, use the convergence ${u}_n\rightharpoonup \widetilde{u}$ in $X$ and the continuity of $j$ to infer that
	\begin{equation}\label{Econv}
	u_n\to \widetilde{u}\quad{\rm in}\quad X.	
	\end{equation}
We now combine \eqref{Ereg} and \eqref{Econv}
to see that $d(u_n,K)\to 0$ which concludes the proof of this step.

\medskip\noindent {\it Step
	{iii)  We  prove that
		any weakly convergent subsequence of the sequence $\{u_n\}$   converges to the solution $u$ of inequality $(\ref{1})$.}} This step is a direct consequence of Steps i), ii) and Theorem \ref{tm}.

	\medskip\noindent {\it Step
		{iv)  We  prove that  the whole sequence $\{u_n\}$ converges to the solution $u$ of inequality $(\ref{1})$. }} To this end, we argue by contradiction. Assume that  the convergence (\ref{c0}) does not hold.  Then, there exists $\delta_0>0$ such that for all $k\in\mathbb{N}$ there exists $u_{n_k}\in X$ such that
		\begin{equation}\label{m6}
		\|u_{n_k}-u\|_X\ge \delta_0.
		\end{equation}
		Note that the sequence $\{u_{n_k}\}$ is a subsequence of the sequence $\{u_n\}$ and, therefore, Step i) and Lemma \ref{lm} imply that it is  bounded in $X$. We now use a compactness argument to deduce that there exists a subsequence of the sequence $\{u_{n_k}\}$, again denoted by $\{u_{n_k}\}$, which is weakly convergent in $X$. Then, Step iii) guarantees that $u_{n_k}\to u$ as $k\to \infty$. We now pass to the limit when  $k\to \infty$ in \eqref{m6} and find that $\delta_0\le 0$. This contradicts inequality $\delta_0>0$ and concludes the proof.
	\hfill$\Box$

\medskip
Note that Corollaries \ref{cor1} and \ref{cor2} represent a version of Theorem \ref{t2} b) and c), obtained under assumption  (\ref{A}) and (\ref{jj}) instead to (\ref{Ap}) and (\ref{j}), respectively. Similar arguments can be used  to prove
the result in Theorem \ref{t2} a). 

\medskip\noindent
{\bf Well-posedness results.} Convergence results for the solution of optimization problems and variational inequalities are strongly related to the well-posedness of these problems. References in the field are  \cite{DZ,LP,LP1,LP2,Ty}  and, more recently \cite{S}. Here we restrict ourselves to mention only two classical  well-posedness concepts in the study of the variational inequality (\ref{1}) and,
to this end, we recall the following definitions.

\begin{Definition}\label{Cd41}
{\rm a)} A sequence $\{u_n\}\subset X$ is called an {\it approximating sequence} for inequality $(\ref{1})$ if
there exists a sequence   $0\le\ve_n\to 0$ such that
\begin{eqnarray}
	&&\label{C42n} u_n\in K,\quad (Au_n,v-u_n)_X+j(v)-j(u_n)+\ve_n\|v-u_n\|_X\\ [2mm]
	&&\qquad\qquad\qquad\ge (f,v-u_n)_X\quad\ \forall\, v\in K,\ n\in\mathbb{N}.\nonumber
\end{eqnarray}

\medskip
{\rm b)}  Inequality $(\ref{1})$ is {\it well-posed in the sense of Tykhonov} (or, equivalently, is {\it Tykhonov well-posed}) if it has  a unique solution and any approximating sequence converges in $X$ to $u$.
\end{Definition}

\begin{Definition}\label{Cd41n}

{\rm a)} A sequence $\{u_n\}\subset X$ is called a {\it generalized (or $LP$) approximating sequence} for inequality $(\ref{1})$ if
there  exist two sequences  $\{w_n\}\subset X$ and  $\{\ve_n\}\subset\real_+$ such that $w_n\to 0_X$ in $X$, $\ve_n\to 0$
and, moreover,
\begin{eqnarray}
	&&\hspace{-12mm}
	\label{C42m}
	u_n+w_n\in K,\quad (Au_n,v-u_n)_X+j(v)-j(u_n)+\ve_n\|v-u_n\|_X\\ [2mm]
	&&\hspace{2mm}\qquad\qquad\ge (f,v-{u_n})_X\quad\ \forall\, v\in K,\ n\in\mathbb{N}.\nonumber
\end{eqnarray}

\medskip
{\rm b)} Inequality $(\ref{1})$ is {\it well-posed in the sense of Levitin-Polyak} (or, equivalently, is {\it Levitin-Polyak well-posed}) if it has a unique solution and any $LP$-approximating sequence converges in $X$ to $u$.
\end{Definition}

It is easy to see that if (\ref{C42n})  holds then both conditions (\ref{c2})(a) and (\ref{c2})(b) are satisfied. Therefore, using Definition \ref{Cd41} and Theorem \ref{tm} it is easy to deduce the following result.

\begin{Corollary}\label{cor3}
	Assume $(\ref{K})$, $(\ref{f})$,  $(\ref{A})$  and $(\ref{jj})$.	Then inequality \eqref{1}  is  well-posed in the sense of Tykhonov.
\end{Corollary}

The Levitin-Polyak well-posedness of inequality \eqref{1} is a consequence of Theorem~\ref{tm}, too, as it follows from the following
result.

\begin{Corollary}\label{cor4}
	Assume $(\ref{K})$, $(\ref{f})$,  $(\ref{A})$  and $(\ref{jj})$.	Then inequality \eqref{1}  is Levitin-Polyak well-posed.
\end{Corollary}

\noindent
{\it Proof.} Let $\{u_n\}\subset X$ be a generalized approximating sequence. Then, Definition \ref{Cd41n} (a) shows that $d(u_n,K)\le\|w_n\|_X$ for each $n\in\mathbb{N}$
and, since $w_n\to 0_X$, we deduce that condition
(\ref{c2})(a) is satisfied. Moreover, inequality (\ref{C42n}) shows that condition
(\ref{c2})(b) is satisfied, too. We now use Theorem \ref{tm} to see that $u_n\to u$ in $X$ and conclude the proof by using Definition \ref{Cd41n} (b).
\hfill$\Box$

\medskip
Definitions \ref{Cd41} a) and \ref{Cd41n} a) show that any approximating sequence
is a generalized approximating sequence, too.
Therefore, using Corollary \ref{cor4} we obtain the following implications.
\begin{eqnarray*}
	&&\mbox{$\{u_n\}$ is an approximating sequence}\\ [2mm]
	&& \Longrightarrow\ \
	\mbox{$\{u_n\}$ is a generalized approximating sequence}\\ [2mm]
	&&\Longrightarrow\ \  u_n\to u\quad{\rm  in}\quad X.
\end{eqnarray*}
The following elementary examples show that  the converse of these implications are not valid.

\begin{Example}\label{ex11}
Consider Problem $\cP$ in the particular case
$X=\real$, $K=[0,1]$, $A=I_X$, $j\equiv 0$ and $f=2$. Then	  $(\ref{1})$ becomes: find $u$ such that
\begin{equation}\label{za}
u\in [0,1],\quad u(v-u)\ge f(v-u)\qquad\forall\,v\in[0,1].
\end{equation}
The solution of this inequality is $u=P_Kf=1$. Let the sequence $\{u_n\}\subset \real$ be given by $u_n=1-\frac{1}{n}$ for all $n\in\mathbb{N}$.
Then $u_n\to u$ but $\{u_n\}$ is not a generalized approximating sequence for $(\ref{za})$. Indeed, assume that there exists $0\le\ve_n\to 0$ such that, for all $n\in\mathbb{N}$
\begin{equation}\label{zan}
u_n(v-u_n)+\ve_n|v-u_n|\ge f(v-u_n)\qquad\forall\,v\in[0,1].
\end{equation}
We fix $n\in\mathbb{N}$ and take
$v=1-\frac{1}{2n}$ in \eqref{zan}. Then, using equalities  $u_n=1-\frac{1}{n}$, $f=2$  we deduce that $\ve_n\ge 1+\frac{1}{n}$. This inequality is valid for each  $n\in\mathbb{N}$,
which contradicts the convergence $\ve_n\to 0$.
\end{Example}

\begin{Example}\label{ex12}
	Consider Problem $\cP$ in the particular case
	$X=\real$, $K=[0,1]$, $A=I_X$, $j\equiv 0$ and $f=1$ and note that	the solution of the corresponding inequality $(\ref{1})$ is $u=P_Kf=1$.
	Let $\{u_n\}\subset \real$ be the sequence given by $u_n=1+\frac{1}{n}$ for all $n\in\mathbb{N}$.
	Then, $\{u_n\}$ is not an approximating sequence, since condition $u_n\in K$ for each $n\in\mathbb{N}$ is not satisfied.
	Nevertheless, $\{u_n\}$ is a
	generalized approximating sequence for $(\ref{za})$. Indeed, it is easy to see that conditions in Definition $\ref{Cd41n}$ {\rm a)} hold with $w_n=-\frac{1}{n}$ and  $\ve_n=\frac{1}{n}$, for all $n\in\mathbb{N}$.
\end{Example}

The examples above show that Tykhonov and Levitin-Polyak well-posedness concepts are not optimal,  in the sense that the approximating sequences and the generalized approximate sequences they generate, respectively, do not recover all the sequences of $X$ which converge to the solution $u$ of the variational inequality (\ref{1}).  This remark leads in a natural way to the following question: how to identify a class of sequences,  say the class of $\cT$-approximating sequences, such that the following equivalence holds:
\begin{eqnarray*}
	&&\mbox{$\{u_n\}$ is a $\cT$-approximating sequence}
	\ \ \Longleftrightarrow\ \   u_n\to u\quad{\rm  in}\quad X.
\end{eqnarray*}

A possible answer to this question is provided by the following definition.

\begin{Definition}\label{CT}
	 A sequence $\{u_n\}\subset X$ is called a {\it $\cT$-approximating sequence} for inequality $(\ref{1})$ if
	there exists a sequence   $0\le\ve_n\to 0$ such that
	\begin{eqnarray}
		&&\label{C43} d(u_n,K)\le\ve_n,\quad (Au_n,v-u_n)_X+j(v)-j(u_n)+\ve_n(1+\|v-u_n\|_X)\\ [2mm]
		&&\qquad\qquad\qquad\ge (f,v-u_n)_X\quad\ \forall\, v\in K,\ n\in\mathbb{N}.\nonumber
	\end{eqnarray}
\end{Definition}

Inspired by Definitions \ref{Cd41} b) and \ref{Cd41n} b) we complete  Definition \ref{CT} as follows.

\begin{Definition}\label{CTn}
	Inequality $(\ref{1})$ is {\it $\cT$-well-posed }  if it has  a unique solution and any $\cT$-approximating sequence converges in $X$ to $u$.
\end{Definition}

Adopting these definitions, we are in a position to state the following two theorems, which represent equivalent formulations of Theorem \ref{tm}.

\begin{Theorem}\label{tmn}
	Assume $(\ref{K})$, $(\ref{f})$,  $(\ref{A})$  and $(\ref{jj})$. Then, a sequence $\{u_n\}\subset X$ converges to the solution of inequality
	$(\ref{1})$ if and only if it is a $\cT$-approximating sequence.
\end{Theorem}

\begin{Theorem}\label{tmnp} 	
	Assume $(\ref{K})$, $(\ref{f})$,  $(\ref{A})$  and $(\ref{jj})$.  Then, inequality $(\ref{1})$ has a unique solution if and only if it is $\cT$-well-posed.
\end{Theorem}

Note that Definition \ref{CTn} introduces a new concept of well-posedness for the variational inequality (\ref{1}). It can be extended to the study of various nonlinear problems like hemivariational inequalities, inclusions, minimization problems, various classes of time-dependent and evolutionary inequalities. Details can be found in the recent book \cite{S}. There, the concept of Tykhonov triple, denoted by $\cT$,  was introduced. Moreover, the main properties of Tykhonov triples have been stated and proved, together with various examples and counter examples.
Then, given a metric space $(X,d)$,  Problem $\cP$
and a  Tykhonov triple $\cT$, both defined on $X$, the abstract concept of $\cT$-well-posedness for Problem $\cP$ has been introduced, { based on two main ingredients: the existence of a unique solution to Problem $\cP$, and the convergence to it to a special kind of sequences, the so-called $\cT$-approximating sequences.
Moreover, various applications in Functional Analysis and Contact Mechanics have been provided.

{
We end this section with the remark that well-posedness concepts can be extended  to abstract problems  for which the set of solutions (assumed to be not empty) is not reduced to a singleton.  For such  problems various concepts of generalized well-posedness have been introduced in the literature. They are based on the definition of a family of so-called generalized approximating sequences and the condition that every sequence of this family has a subsequence which converges to some point of the solution set. A recent reference on this topic is \cite{Dey}, where
the Levitin-Polyak well-posedness of the so-called split equilibrium problems is studied. Additional details can be found in the book \cite{S}.}

}

\section{A frictionless contact problem}\label{s5}
\setcounter{equation}0

The abstract results  in Sections \ref{s3} and \ref{s4} are useful
in the study of various mathematical models which describe the equilibrium of elastic bodies in  contact with an obstacle, the so-called foundation. In this section
we introduce and study an example of such model and, to this end, we need some  notations and preliminaries.

Let $d\in\{2,3\}$. We denote by $\mathbb{S}^d$ the space of second order symmetric tensors on $\mathbb{R}^d$ and use  the notation $``\cdot"$, $\|\cdot\|$, $\bzero$ for the inner product, the norm and the zero element of the spaces
$\mathbb{R}^d$ and $\mathbb{S}^d$, respectively.
Let $\Omega\subset\mathbb{R}^d$ be a domain with smooth boundary $\Gamma$ divided into three
measurable disjoint parts $\Gamma_1$, $\Gamma_2$ and $\Gamma_3$ such that ${ meas}\,(\Gamma_1)>0$.
A generic point in $\Omega\cup\Gamma$ will be denoted by $\bx=(x_i)$.
We use the
standard notation for Sobolev and Lebesgue spaces associated to
$\Omega$ and $\Gamma$. In particular, we use the spaces  $L^2(\Omega)^d$, $L^2(\Gamma_2)^d$
and $H^1(\Omega)^d$, endowed with their canonical inner products and associated norms.
Moreover, for an element $\bv\in H^1(\Omega)^d$ we still  write $\bv$ for the trace of
$\bv$ to $\Gamma$. In addition, we consider the
space
\begin{eqnarray*}
	&&V=\{\,\bv\in H^1(\Omega)^d\ :\  \bv =\bzero\ \ {\rm a.e.\ on\ \ }\Gamma_1\,\},
\end{eqnarray*}
which is a real Hilbert space
endowed with the canonical inner product
\begin{equation*}
	(\bu,\bv)_V= \int_{\Omega}
	\bvarepsilon(\bu)\cdot\bvarepsilon(\bv)\,dx
\end{equation*}
and the associated norm
$\|\cdot\|_V$. Here and below $\bvarepsilon$
represents the deformation operator, i.e.,
\[
\bvarepsilon(\bu)=(\varepsilon_{ij}(\bu)),\quad
\varepsilon_{ij}(\bu)=\frac{1}{2}\,(u_{i,j}+u_{j,i}),
\]
where an index that follows a comma denotes the
partial derivative with respect to the corresponding component of $\bx$, e.g.,\ $u_{i,j}=\frac{\partial u_i}{\partial x_j}$.
The completeness of the space $V$ follows from the assumption
${ meas}\,(\Gamma_1)>0$ which allows us to use Korn's inequality.
We denote by $\bzero_V$ the zero element of $V$ and we recall that, for an element $\bv\in V$,  its  normal and tangential components on $\Gamma$
are given by
\[\mbox{$v_\nu=\bv\cdot\bnu$\ \quad  and\ \quad  $\bv_\tau=\bv-v_\nu\bnu$,}\] respectively.   Here and below $\bnu$ denote the unitary outward normal to $\Gamma$.
We also  recall the trace inequality
\begin{equation}\label{trace}
	\|\bv\|_{L^2(\Gamma)^d}\leq d_0\|\bv\|_{V}\qquad \forall\,
	\bv\in V
\end{equation}
in which $d_0$ represents a positive constant.

\smallskip
For the inequality problem we consider in this section we use the data ${\cal F}$, $F$, $\fb_0$, $\fb_2$ and $k$ which satisfy
the following conditions.

\begin{eqnarray}
	&&\label{Fc}\left\{\begin{array}{ll}
		{\rm (a)}\ {\cal F}\colon
		\mathbb{S}^d\to \mathbb{S}^d. \\ [1mm]
		{\rm (b)\  There\ exists}\ M_{\cal F}>0\ {\rm such\ that}\\
		{}\qquad \|{\cal F}\bvarepsilon_1-{\cal F}\bvarepsilon_2\|
		\le M_{\cal F} \|\bvarepsilon_1-\bvarepsilon_2\|\quad\mbox{for all} \ \ \bvarepsilon_1,\bvarepsilon_2
		\in \mathbb{S}^d.
		\\ [1mm]
		{\rm (c)\ There\ exists}\ m_{\cal F}>0\ {\rm such\ that}\\
		{}\qquad ({\cal F}\bvarepsilon_1-{\cal F}\bvarepsilon_2)
		\cdot(\bvarepsilon_1-\bvarepsilon_2)\ge m_{\cal F}\,
		\|\bvarepsilon_1-\bvarepsilon_2\|^2\quad \mbox{for all} \ \ \bvarepsilon_1,
		\bvarepsilon_2 \in \mathbb{S}^d.
	\end{array}\right.
\\ [2mm]
	&&\label{F}F\in L^\infty(\Gamma_3),\qquad F(\bx)\ge 0\ \ {\rm a.e.}\ \bx\in \Gamma_3.\\ [2mm]
	&&\label{f0}\fb_0\in L^2(\Omega)^d,\qquad\fb_2\in L^2(\Gamma_2)^d.\\ [2mm]
	&&k>0. \label{k}
\end{eqnarray}
Moreover, we use $K$ for the set defined by
\begin{equation}
	\label{KK}K=\{\,\bv\in V\ :\ v_\nu \le k\ \  \hbox{a.e. on}\
	\Gamma_3\,\}
\end{equation}
and  $r^+$ for the positive part of $r\in\mathbb{R}$, that is $r^+=\max\,\{r,0\}$.

\medskip
Then, the inequality problem we consider  is the following.

\medskip\medskip\noindent
\medskip\noindent{\bf Problem}  ${\cal P}^c$. {\it Find  $\bu$
	such that}
\begin{eqnarray*}\label{51}
	&&\bu\in K,\quad \int_{\Omega}{\cal F}\bvarepsilon(\bu)\cdot(\bvarepsilon(\bv)-\bvarepsilon(\bu))\,dx+\int_{\Gamma_3}Fu_\nu^+(v_\nu-u_\nu)\,da\\[2mm]
	&&\qquad\quad\ge
	\int_{\Omega}\fb_0\cdot(\bv-\bu)\,dx
	+\int_{\Gamma_2}\fb_2\cdot(\bv-\bu)\,da\quad\forall\,\bv\in K.\nonumber
\end{eqnarray*}

\medskip
Following the arguments in \cite{SofMat,SMBOOK}, it can be shown that Problem $\cP^c$ represents the variational formulation of a mathematical model that describes the equilibrium of an elastic body  $\Omega$ which is acted upon by  external forces,
is fixed on $\Gamma_1$, and is
in frictionless contact on $\Gamma_3$. The contact takes place with a rigid foundation  covered by a layer of rigid-plastic material of thickness $k$.
Here ${\cal F}$ is the elasticity operator,   $\fb_0$ and $\fb_2$ denote the density of  applied body forces and tractions which act on the body and the surface $\Gamma_2$, respectively and $F$ is a given function which describes the yield limit of the rigid-plastic material.

Next, consider the sequences  $\{F_n\}$, $\{\mu_n\}$, $\{\fb_{0n}\}$, $\{\fb_{2n}\}$, $\{k_n\}$ such that, for each $n\in\mathbb{N}$, the following hold.
\begin{eqnarray}	
&&\label{Fn} F_n\in L^\infty(\Gamma_3),\qquad F_n(\bx)\ge 0\ \ {\rm a.e.}\ \bx\in \Gamma_3.\\ [2mm]
&&\label{mu}\mu_n\in L^\infty(\Gamma_3),\qquad \mu_n(\bx)\ge 0\ \ {\rm a.e.}\ \bx\in \Gamma_3.\\ [2mm]
&&\label{f0n}\fb_{0n}\in L^2(\Omega)^d,\qquad\fb_{2n}\in L^2(\Gamma_2)^d.\\ [2mm]
&&k_n\ge k. \label{kn}\\ [2mm]
&&d_0^2\|\mu\|_{L^\infty(\Gamma_3)} \|F_n\|_{L^\infty(\Gamma_3)}<m_{\cal F}. \label{sm}
\end{eqnarray}
Recall that in (\ref{sm}) and below $d_0$ and $m_{\cal F}$ represent the constants introduced in \eqref{trace} and \eqref{Fc}, respectively.
Finally, for each $n\in\mathbb{N}$
define the set
\begin{equation}
	\label{KKn}K_n=\{\,\bv\in V\ :\ v_\nu \le k_n\ \  \hbox{a.e. on}\
	\Gamma_3\,\}.
\end{equation}
and  consider the following perturbation of Problem $\cP^c$

\medskip\medskip\noindent{\bf Problem}  ${\cal P}^c_n$. {\it Find $\bu_n$ such that}				
\begin{eqnarray*}\label{51n}
	&&\bu_n\in K_n,\quad \int_{\Omega}{\cal F}\bvarepsilon(\bu_n)\cdot(\bvarepsilon(\bv)-\bvarepsilon(\bu_n))\,dx+\int_{\Gamma_3}Fu_{n\nu}^+(v_\nu-u_{n\nu})\,da\\[3mm]
	&&\qquad+\int_{\Gamma_3}\mu_n\,F_nu_{n\nu}^+(\|{\bv}_\tau\|-\|{\bu_n}_\tau\|)\,da\nonumber\\ [3mm]
	&&\qquad\qquad\ge
	\int_{\Omega}\fb_{0n}\cdot(\bv-\bu)\,dx
	+\int_{\Gamma_2}\fb_{2n}\cdot(\bv-\bu_n)\,da\quad\forall\,\bv\in K_n.\nonumber
\end{eqnarray*}

\medskip
Problem $\cP_n^c$ represents the variational formulation of a mathematical model of contact, similar to that associated to Problem $\cP^c$. Nevertheless, two differences exist between the corresponding models. The first one arises in the fact that for the model in Problem $\cP_n^c$ the contact is assumed to be frictional and is described with the classical Coulomb law of dry friction, governed by the friction coefficient $\mu_n$. The second  difference arises in the fact that in the statement of Problem $\cP^c_n$  we use the data $\fb_{0n}$, $\fb_{2n}$ and $k_n$ which represent a perturbation of the data
$\fb_{0}$, $\fb_{2}$ and $k$, respectively, used in the statement of Problem $\cP^c$.

\medskip

Our main result in this section,  is the following.

\begin{Theorem}\label{t5}
	Assume $(\ref{Fc})$--$(\ref{k})$, $(\ref{Fn})$--$(\ref{sm})$. Then
	Problem $\cP^c$  has a unique solution  and, for each $n\in\mathbb{N}$, Problem $\cP^c_n$ has a unique solution.
	Moreover, if
	\begin{equation}\label{co1c}
		\left\{\begin{array}{ll}
		k_n\to k,\quad
		\mu_n\to 0\ \ {\rm in}\ \ L^\infty(\Gamma_3),\quad
		F_n\to F\ \ {\rm in}\ \ L^\infty(\Gamma_3),\\ [5mm]
		\fb_{0n}\to \fb_{0}\ \ {\rm in}\ \ L^2(\Omega)^d,\quad \fb_{2n}\to \fb_2\ \ {\rm in}\ \ L^2(\Gamma_2)^d\quad{\rm as}\ \ n\to \infty,
		\end{array}\right.
	\end{equation}
 then the solution of Problem $\cP^c_n$ converges to the solution of Problem $\cP^c$, i.e.,
	\begin{equation}\label{co2c}
		\bu_n\to \bu\quad {\rm in}\ \ V\quad{\rm as}\quad n\to \infty.
	\end{equation}
\end{Theorem}

\medskip\noindent{\it Proof.}
We start with some additional notation. First,  we
consider the operator $A:V\to V$, the functions $j,\, j_n:V\to\R$, $\varphi_n:V\times V\to\R$ and the elements $\fb,\, \fb_n\in V$ defined as follows:

\begin{eqnarray}
	&&
	\label{8b1}(A\bu,\bv)_V =\int_{\Omega}\cF\bvarepsilon(\bu)\cdot\bvarepsilon(\bv)\,dx,\\ [2mm]
	&&
	\label{8b3}
	j(\bv)=\int_{\Gamma_3}Fv_\nu^+\,da,\quad
	j_n(\bv)=\int_{\Gamma_3}F_nv_\nu^+\,da,\\ [2mm]
	&&
	\label{8b3b}
	\varphi_n(\bu,\bv)=\int_{\Gamma_3}\mu_n F_nu_\nu^+\|{\bv}_\tau\|\,da,\\ [2mm]
	&&
	\label{8b4}\left\{\begin{array}{ll}(\fb,\bv)_V=\displaystyle\int_{\Omega}\fb_0\cdot\bv\,dx
	+\int_{\Gamma_2}\fb_2\cdot\bv\,da,\\[4mm](\fb_n,\bv)_V=\displaystyle\int_{\Omega}\fb_{0n}\cdot\bv\,dx
	+\int_{\Gamma_2}\fb_{2n}\cdot\bv\,da
	\end{array}\right.
\end{eqnarray}
for all $\bu,\bv\in V$ and $n\in\mathbb{N}$. Then, it is easy to see that

\begin{equation}\label{e1c}
	\left\{\begin{array}{l}
		\mbox{$\bu$ is a solution of Problem $\cP^c$ if and only if}\\ [4mm]
		\bu\in K, \quad (A\bu,\bv-\bu)_V+j(\bv)-j(\bu)\ge (\fb,\bv-\bu)_V\quad \forall\, \bv\in K.
	\end{array}\right.
\end{equation}

Moreover, for each $n\in\mathbb{N}$, the following equivalence holds:
\begin{equation}\label{e2c}
	\left\{\begin{array}{l}
		\mbox{$\bu_n$ is a solution of Problem $\cP^c_n$ if and only if\qquad}\\ [4mm]
		\bu_n\in K_n, \quad (A\bu_n,\bv-\bu_n)_V+j_n(\bv)-
		j_n(\bu_n)\\ [2mm]
		\qquad\qquad\qquad+\varphi_n(\bu_n,\bv)-\varphi_n(\bu_n,\bu_n)\ge (\fb_n,\bv-\bu_n)_V\quad \forall\, \bv\in K_n.
	\end{array}\right.
\end{equation}

\medskip
Equivalence \eqref{e1c} suggests us to use the abstract results in Sections \ref{s2} and \ref{s3} with $X=V$, $K$ defined by (\ref{KK}), $A$ defined by (\ref{8b1}), $j$ defined by (\ref{8b3}) and $\fb$ given by (\ref{8b4}). It is easy to see that in this case conditions $(\ref{K})$, $(\ref{f})$, $(\ref{A})$ and $(\ref{jj})$  are satisfied.
For instance, using assumption (\ref{Fc})  we see that
\begin{eqnarray*}
	(A\bu - A\bv,\bu -\bv)_{V} \geq m_{\cal F} \|\bu -\bv\|^{2}_{V},\qquad \|A\bu-A\bv\|_V\le M_{\cal F}\, {\|\bu-\bv\|_V}
\end{eqnarray*}
for all $\bu,\, \bv\in V$. Therefore, conditions (\ref{A1}) and  (\ref{A2}) hold with $m=m_{\cal F}$ and $M=M_{\cal F}$, respectively which shows that $A$ satisfies condition (\ref{A}).
Condition (\ref{jj}) is  also satisfied
since $j$ is a continuous seminorm on the space $V$.

Therefore, we are in a position to apply Theorem \ref{t1}  in order to deduce the existence of a unique solution of the variational inequality in (\ref{e1c}).
The unique solvability of the variational inequality in (\ref{e2c}) follows from a standard argument of quasivariational inequalities. The proof can be found in \cite{SofMat}, for instance and, therefore, we skip it. Note that  here, besides the regularities (\ref{Fn})--(\ref{f0n}) and condition $k_n>0$, we need the smallness assumption (\ref{sm}).

We now move to the proof of the convergence \eqref{co2c}. Let $n\in\mathbb{N}$ and $\bv\in K$. We write
\begin{eqnarray}
&&\label{w1}(A\bu_n,\bv-\bu_n)_V+j(\bv)-j(\bu_n)- (\fb,\bv-\bu_n)_V
\\ [2mm]
&&=(A\bu_n,\bv-\bu_n)_V+j_n(\bv)-j_n(\bu_n)
+\varphi_n(\bu_n,\bv)-\varphi_n(\bu_n,\bu_n)-(\fb_n,\bv-\bu_n)_V
\nonumber\\ [2mm]
&&+\big[j(\bv)-j(\bu_n)-j_n(\bv)+j_n(\bu_n)\big]
+\big[\varphi_n(\bu_n,\bu_n)-\varphi_n(\bu_n,\bv)\big]+(\fb_n-\fb,\bv-\bu_n)_V.\nonumber
\end{eqnarray}
Then, we use assumption $k_n\ge k$ to see that
$K\subset K_n$ and, therefore, we are allowed to test in  (\ref{e2c}) with $\bv\in K$. We obtain
\begin{equation}\label{w2}
(A\bu_n,\bv-\bu_n)_V+j_n(\bv)-j_n(\bu_n)
+\varphi_n(\bu_n,\bv)-\varphi_n(\bu_n,\bu_n)-(\fb_n,\bv-\bu_n)_V\ge 0
\end{equation}
and, combining \eqref{w1} and \eqref{w2}
we deduce that
\begin{eqnarray}
	&&\label{w3}(A\bu_n,\bv-\bu_n)_V+j(\bv)-j(\bu_n)+ \big[j(\bu_n)-j(\bv)+j_n(\bv)-j_n(\bu_n)\big]\\ [2mm]
	&&\qquad+\big[\varphi_n(\bu_n,\bv)-\varphi_n(\bu_n,\bu_n)\big]+(\fb-\fb_n,\bv-\bu_n)_V\ge
	(\fb,\bv-\bu_n)_V.\nonumber
\end{eqnarray}

Next, using definitions (\ref{8b3})--(\ref{8b4}) and standard embedding and trace arguments we find that
\begin{eqnarray}
&&\label{w4}j(\bu_n)-j(\bv)+j_n(\bv)-j_n(\bu_n)\le c_0\|F_n-F\|_{L^\infty(\Gamma_3)}\|\bv-\bu_n\|_V,\\[3mm]
&&\label{w5}(\fb-\fb_n,\bv-\bu_n)_V\\
&& \qquad\le c_0\big(\|\fb_{0n}-\fb_0\|_{L^2(\Omega)^d}+
\|\fb_{2n}-\fb_2\|_{L^2(\Gamma_2)^d}\big)\|\bv-\bu_n\|_V,\nonumber\\ [3mm]
&&\label{w6}\varphi_n(\bu_n,\bv)-\varphi_n(\bu_n,\bu_n)\le
c_0\|\mu_n\|_{L^\infty(\Gamma_3)}\|F_n\|_{L^\infty(\Gamma_3)}\|\bv-\bu_n\|_V,
\end{eqnarray}
where $c_0$ represents a positive constant which does not depend on $n$. We now substitute
\eqref{w4}--\eqref{w6} in \eqref{w3} and use notation
\begin{eqnarray}
&&\label{w7}\ve_n={\rm max}\,\Big\{c_0\|F_n-F\|_{L^\infty(\Gamma_3)},c_0\big(\|\fb_{0n}-\fb_0\|_{L^2(\Omega)^d}+
\|\fb_{2n}-\fb_2\|_{L^2(\Gamma_2)^d}\big),\nonumber \\ [2mm]
&&\qquad\qquad\qquad c_0\|\mu_n\|_{L^\infty(\Gamma_3)}\|F_n\|_{L^\infty(\Gamma_3)}\Big\}
\end{eqnarray}
to deduce that
\begin{equation}
	\label{w8}(A\bu_n,\bv-\bu_n)_V+j(\bv)-j(\bu_n)+ \ve_n\|\bv-\bu_n\|_V\ge
	(\fb,\bv-\bu_n)_V.
\end{equation}
Note that definition \eqref{w7} and assumptions
\eqref{co1c} imply that $\ve_n\to 0$. Therefore,
inequality \eqref{w8} shows that condition (\ref{c2})(b) is satisfied.  We are now in a position to use Lemma \ref{lm} to see that
\begin{equation}\label{w9}
\mbox{there exists $D>0$ such that}\ \|\bu_n\|_V\le D\quad\forall\, n\in\mathbb{N}.
\end{equation}

On the other hand, using definitions \eqref{KK} and \eqref{KKn} it is easy to see that $\frac{k}{k_n}\bu_n\in K$ and, therefore

\begin{equation}\label{w10}d(\bu_n,K)\le \Big\|\bu_n-\frac{k}{k_n}\bu_n\Big\|_V=\Big|1-\frac{k}{k_n}\Big|\|\bu_n\|_V \quad\forall\, n\in\mathbb{N}.
\end{equation}
We now use \eqref{w10}, \eqref{w9} and assumption $k_n\to k$ to see that
$d(\bu_n,K)\to 0$ which shows that condition (\ref{c2})(a) is satisfied, too.

It follows from above that we are in a position to use Theorem \ref{tm} to deduce the convergence (\ref{co2c}).
These results combined with equivalences (\ref{e1c}) and (\ref{e2c}) allows us to conclude the proof of the  theorem.
\hfill$\Box$

\medskip
We end this section with the following physical interpretation of Theorem \ref{t5}.
First, the existence and uniqueness part in the theorem proves the unique weak solvability of the contact problems considered: the frictionless contact with a rigid foundation covered by a layer of rigid-plastic material of thickness $k$ and the frictional contact with a rigid foundation covered by a layer of rigid-plastic material of thickness $k_n$.
Second,
the weak solution of the frictionless contact problem with a rigid foundation  covered by a layer of rigid-plastic material depends continuously on the density of body forces and
surface tractions as well as on the yield limit and the thickness of the layer. In addition, it can be approached by the solution of the corresponding frictional problem with a small coefficient of friction.

\section{A heat transfer problem}\label{s6}
\setcounter{equation}0

In this section we apply the abstract results in Sections \ref{s2}--\ref{s4} in the study of a mathematical model which describes a heat transfer phenomenon.  The problem we consider represents a version of the problems already considered in \cite{GMOT,ST3} and, for this reason, we skip the details. Its classical formulation  is the following.

\medskip\noindent{\bf Problem $\cC^t$}. {\it
	Find a temperature field $u:\Omega\to\R$ such that}
\begin{eqnarray}
	&&\label{d1} -\Delta u=g\qquad\ \ {\rm a.e.\ in\ }\Omega,\\ [3mm]
	&&\label{d2}u=0\hspace{18mm}{\rm a.e.\ on\ }\Gamma_1,\\ [3mm]
	&&\label{d3}-\frac{\partial u}{\partial\nu}=q\hspace{11mm}{\rm a.e.\ on\ }\Gamma_2.\\ [3mm]
	&&\label{d4}u=b\hspace{18mm}{\rm a.e.\ on\ }\Gamma_3.
	\end{eqnarray}

\medskip
Here, as in Section \ref{s5},  $\Omega$ is a bounded domain in $\real^d$ ($d=2,3$ in applications) with smooth boundary
$\Gamma$, divided into three
measurable disjoint parts $\Gamma_1$, $\Gamma_2$ and $\Gamma_3$ such that ${ meas}\,(\Gamma_1)>0$.   We denote by $\bnu$ and outer normal unit to $\Gamma$
and recall that
in (\ref{d1})--(\ref{d4}) we do not mention the dependence of the different functions on the spatial variable $\bx\in\Omega\cup\Gamma$.
The functions $g$, $q$ and $b$ are given and will be described below. Here we mention that  $g$ represents the internal energy,
 $q$ is a prescribed  heat flux and $b$ denotes a prescribed temperature. Moreover, $\frac{\partial u}{\partial\nu}$ denotes the normal derivative of $u$ on the boundary $\Gamma$.

 \medskip
 Now, let $\{\lambda_n\}\subset\real$ be a sequence such that $\lambda_n>0$ for each $n\in\mathbb{N}$. Then for each $n\in\mathbb{N}$ we consider the following boundary problem.

 \medskip\noindent{\bf Problem $\cC^t_n$}. {\it
 	Find a temperature field $u_n:\Omega\to\R$ such that}
 \begin{eqnarray}
 	&&\label{d1n} -\Delta u_n=g\hspace{23mm}{\rm a.e.\ in\ }\Omega,\\ [2mm]
 	&&\label{d2n}u_n=0\hspace{30mm}{\rm a.e.\ on\ }\Gamma_1,\\ [2mm]
 	&&\label{d3n}-\frac{\partial u_n}{\partial\nu}=q\hspace{24mm}{\rm a.e.\ on\ }\Gamma_2,\\ [2mm]
 	&&\label{d4n}-\frac{\partial u_n}{\partial\nu}=\frac{1}{\lambda_n}(u_n-b)\hspace{6mm}{\rm a.e.\ on\ }\Gamma_3.
 \end{eqnarray}

Note that Problem ${\cC}^t_n$ is obtained from Problem $\cC^t$ by replacing the Dirichlet boundary condition (\ref{d4}) with the Neumann boundary condition (\ref{d4n}). Here $\lambda_n$ is a positive parameter, and its inverse $h_n=\frac{1}{\lambda_n}$ represents the heat transfer coefficient on the boundary $\Gamma_3$. In contrast to Problem  $\cC^t$ (in which the temperature is prescribed on $\Gamma_3$), in Problem $\cC^t_n$ this condition is replaced by a condition on the flux of the temperature, governed by a positive  heat transfer coefficient.

\medskip
For the variational analysis of Problem $\cC^t$  we consider the space
\[V=\ \{\,v\in H^1(\Omega)\ :\ v=0\quad{\rm a.e.\ on}\quad\Gamma_1\}.\]
Note that, here and below in this section, we still write $v$ for trace of the element $v$ to $\Gamma$. Denote in what follows by $(\cdot,\cdot)_V$ the inner product of the space $H^1(\Omega)$ restricted to $V$ and by $\|\cdot\|_V$ the associated norm.
Since $meas \,\Gamma_1>0$, it is well known that $(V,(\cdot,\cdot)_V)$ is a real Hilbert space.
Next, we assume that
\begin{eqnarray}
	&&\label{dw}
	g\in L^2(\Omega),\quad q\in L^2(\Gamma_2),\quad b\in H^{\frac{1}{2}}(\Gamma_3),\\ [2mm]
	&&\label{dwz} \mbox{there exists $v_0\in V$ such that  $v_0=b$\ \ a.e.\ on\ \ $\Gamma_3$,}
\end{eqnarray}
and, finally, we introduce the set
\begin{eqnarray}
	&&\label{d6} K=\ \{\,v\in V\ :\  v=b\ \ {\rm a.e.\ on}\ \ \Gamma_3\,\}.
\end{eqnarray}
Note that assumption (\ref{dwz}) represents a compatibilty assumption on the data $b$ which guarantees that the set $K$ is not empty.
Then, it is easy to see that the variational formulation of problems $\cC^t$ and $\cC^t_n$, obtained by standard arguments, are as follows.

\medskip\medskip\noindent{\bf Problem}  ${\cal P}^t$. {\it Find $u$ such that}
\begin{equation*}\label{1d}u\in K,\quad\int_{\Omega}\nabla u\cdot(\nabla v-\nabla u)\,dx+\int_{\Gamma_2}q(v-u)\,da=\int_\Omega g(v-u) \,dx \quad\forall\,v\in K.
\end{equation*}

\medskip\medskip\noindent{\bf Problem}  ${\cal P}^t_n$. {\it Find $u_n$ such that}
\begin{eqnarray*}
	&&\label{1dx}u_n\in V,\qquad\int_{\Omega}\nabla u_n\cdot(\nabla v-\nabla u_n)\,dx+\int_{\Gamma_2}q(v-u_n)\,da\\[2mm]
	&&\qquad\qquad\ \ +\frac{1}{\lambda_n}
	\int_{\Gamma_3}(u_n-b)(v-u_n)\,da \ge\int_\Omega g(v-u_n) \,dx \quad\forall\,v\in V.\nonumber
\end{eqnarray*}

\medskip  Our main result in this section is the following.

\begin{Theorem}\label{t7}
	Assume $(\ref{dw})$ and $(\ref{dwz})$. Then, Problem $\cP^t$  has a unique solution  and, for each $n\in\mathbb{N}$, Problem $\cP^t_n$ has a unique solution. Moreover, if  $\lambda_n\to 0$, then
	the solution of Problem $\cP^t_n$ converges to the solution of Problem $\cP^t$, i.e.,
	\begin{equation}\label{cot2}
		u_n\to u\quad {\rm in}\ \ V\qquad\quad{\rm as}\quad n\to \infty.
	\end{equation}

\end{Theorem}

\noindent{\it Proof.} We
consider the operators $A:V\to V$, $G:V\to V$ and the element $f\in V$  defined as follows:
\begin{eqnarray}
	&&\label{d8} \ (Au,v)_V=\int_\Omega \nabla u\cdot\nabla v\,dx\qquad\forall\,u, v\in V,\\ [2mm]
	&&\label{d9}  (Gu,v)_V=\int_{\Gamma_3}(u-b)v\,da\qquad\forall\,u, v\in V,\\ [2mm]
	&&\label{dd9} \ (f,v)_V=\int_{\Omega}gv\,dx- \int_{\Gamma_2}qv\,da\qquad\forall\, v\in V. 
\end{eqnarray}
\medskip
Then, since the set $\{\,v-v_0\ :\ v\in K\,\}$ is a linear subspace on $V$, it is easy to see that
\begin{equation}\label{e1}
	\left\{\begin{array}{l}
		\mbox{$u$ is a solution of Problem $\cP^t$ if and only if}\\ [3mm]
		u\in K, \quad (Au,v-u)_V\ge (f,v-u)_{V}\quad\ \forall\, v\in K.
	\end{array}\right.
\end{equation}
Moreover, for each $n\in\mathbb{N}$, \begin{equation}\label{e2}
	\left\{\begin{array}{l}
		\mbox{$u_n$ is a solution of Problem $\cP^t_n$ if and only if\qquad}\\ [3mm]
		u_n\in V, \quad (Au_n,v-u_n)_V+\frac{1}{\lambda_n}(Gu_n,v-u_n)_V\\ [2mm]
		\qquad\qquad\qquad\ge (f,v-u_n)_{V}\quad\ \forall\, v\in V.
	\end{array}\right.
\end{equation}

 We use the abstract results in Sections \ref{s2} and \ref{s3} with $X=V$, $K$ defined by (\ref{d6}), $A$ defined by (\ref{d8}), $G$ defined by (\ref{d9}), $f$ given by \eqref{dd9}, 
and $j\equiv 0$. It is easy to see that in this case conditions $(\ref{K})$, $(\ref{f})$, $(\ref{la})$, $(\ref{G})$, $(\ref{A})$ and $(\ref{jj})$ are satisfied. Therefore, we are in a position to apply  Theorem \ref{t1} in order to deduce the existence of a unique solution of the variational inequalities in (\ref{e1}) and (\ref{e2}), respectively.  Moreover, using Corollary \ref{cor2} we deduce the convergence (\ref{cot2}).
These results combined with (\ref{e1}) and (\ref{e2}) allow us to conclude the proof.
\hfill$\Box$

\medskip
We end this section with the following physical interpretation of Theorem \ref{t7}.
First, the solutions of Problems $\cP^t$ and $\cP^t_n$ represent weak solutions of the heat transfer Problems $\cC^t$ and $\cC^t_n$, respectively. Therefore,
Theorem \ref{t7} provides the unique weak solvability of these problems. Second, the convergence result \eqref{cot2} shows that
the weak solution of Problem $\cC^t$ with prescribed temperature on $\Gamma_3$ can be approached by the solution of Problem $\cC_n^t$ with heat transfer on $\Gamma_3$, for a large heat transfer coefficient.

\section*{Acknowledgement}

\indent This project has received funding from the European Union’s Horizon 2020
Research and Innovation Programme under the Marie Sklodowska-Curie
Grant Agreement No 823731 CONMECH. The second author was also supported by the Ministry of Science and Higher Education of Republic of Poland under Grant No 440328/PnH2/2019, and in part from National Science Centre, Poland under project OPUS no. 2021/41/B/ST1/01636.

\end{document}